\documentclass[12pt]{article}

\usepackage[T1]{fontenc}
\usepackage[active]{srcltx}
\usepackage{lmodern}
\usepackage{amsmath,amstext,amssymb,graphicx,graphics,color}
\usepackage{tikz}
\usepackage{pgfplots}
\usepackage{theorem}
\usepackage{cite}               % to cite in increasing order
\usepackage{hyperref}           % faire un pdf avec des liens dynamiques
\usepackage{bbm}                % \mathbbm{1}
\usepackage{fancybox}           % pour \ovalbox
\usepackage[section]{placeins}  % \FloatBarrier
\usepackage{subcaption}
\theorembodyfont{\normalfont\slshape} % on remet la config normal
\newtheorem{definition}{Definition}
\newtheorem{theorem}{Theorem}
\newtheorem{proposition}{Proposition}[section]

\newtheorem{lemma}[proposition]{Lemma}
\theoremstyle{break} % les theorem vont à la ligne dorénavant ("plain" otherwise)

{{\par\noindent \bf Proof. \nobreak}}%
{\nobreak \removelastskip \nobreak \hfill $\Box$ \medbreak}
{{\par\noindent \bf Proof \nobreak}}%
{\nobreak \removelastskip \nobreak \hfill $\Box$ \medbreak}
{{\par\noindent \bf Proof lemma. \nobreak}}%
{\nobreak \removelastskip \nobreak \bf End proof lemma. \medbreak}
\newenvironment{remark}{\par \medskip \noindent {\bf Remark. }\nobreak}{\par \medskip}
\usepackage{fancyhdr}
 \fancypagestyle{plain}{
  \fancyhead{}
  
}
\def\paragraph#1{{\bf #1\ }}
\newcommand{\expo}{\mathrm{e}}

\newcommand{\Var}{\mathrm{Var}}
\newcommand{\dd}{\mathrm{d}}

\newcommand{\emp}{\mathrm{emp}}
\newcommand{\VV}{\mathrm{V}}
\newcommand{\overbar}[1]{\mkern 1.5mu\overline{\mkern-1.5mu#1\mkern-1.5mu}\mkern 1.5mu}
\usepackage[utf8]{inputenc}
\usepackage{unicode_character}

\def\Proof{\noindent{\bf Proof}\quad}
\def\qed{\hfill$\square$\smallskip}

\title{Explicit decay rate for the Gini index in the repeated averaging model}

\author{Fei Cao\footnotemark[1]}

\begin{document}
\maketitle

\footnotetext[1]{Arizona State University - School of Mathematical and Statistical Sciences, 900 S Palm Walk, Tempe, AZ 85287-1804, USA}

\tableofcontents

\begin{abstract}
We investigate the repeated averaging model for money exchanges: two agents picked uniformly at random share half of their wealth to each other. It is intuitively convincing that a Dirac distribution of wealth (centered at the initial average wealth) will be the long time equilibrium for this dynamics. In other words, the Gini index should converge to zero. To better understand this dynamics, we investigate its limit as the number of agents goes to infinity by proving the so-called propagation of chaos, which links the stochastic agent-based dynamics to a (limiting) nonlinear partial differential equation (PDE). This deterministic description has a flavor of the classical Boltzmann equation arising from statistical mechanics of dilute gases. We prove its convergence toward its Dirac equilibrium distribution by showing that the associated Gini index of the wealth distribution converges to zero with an explicit rate.
\end{abstract}


\section{Introduction}

Econophysics is an emerging branch of statistical physics that apply concepts and techniques of traditional physics to economics and finance \cite{savoiu_econophysics:_2013,chatterjee_econophysics_2007,dragulescu_statistical_2000}. It has attracted considerable attention in recent years raising challenges on how various economical phenomena could be explained by universal laws in statistical physics, and we refer to \cite{chakraborti_econophysics_2011,chakraborti_econophysics_2011-1,pereira_econophysics_2017,kutner_econophysics_2019}  for a general review.

% Part 1. emerging brach... -> review '11 (more up to date?)
% Part 2. present other type of models
The primary motivation for study models arising from econophysics is at least two-fold: from the perspective of a policy maker, it is important to deal with the raise of income inequality \cite{dabla-norris_causes_2015,de_haan_finance_2017} in order to establish a more egalitarian society. From a mathematical point of view, we have to understand the fundamental mechanisms, such as money exchange resulting from individuals, which are usually agent-based models. Given an agent-based model, one is expected to identify the limit dynamics as the number of individuals tends to infinity and then its corresponding equilibrium when run the model for a sufficiently long time (if there is one), and this guiding approach is carried out in numerous works across different fields among literatures of applied mathematics, see for instance \cite{naldi_mathematical_2010,barbaro_phase_2014,carlen_kinetic_2013}.

In this work, we consider the so-called repeated averaging model for money exchange in a closed economic system with $N$ agents. The dynamics consists in choosing at random time two individuals and to redistribute equally their combined wealth. To write this dynamics mathematically, we denote by $X_i(t)$ the amount of dollar agent $i$ has at time $t$ for $1\leq i\leq N$. At a random time generated by a Poisson clock with rate $N$, two agents (say $i$ and $j$) update their wealth according to the following rule:
\begin{equation}
\label{eq:repeated_averaging}
\big(X_i,X_j\big) \leadsto \left(\frac{X_i\!+\!X_j}{2}\,,\,\frac{X_i\!+\!X_j}{2}\right),
\end{equation}

Despite of the simplicity of the model, there are actually quite a few manuscripts in the literature which are solely dedicated to it. To the best of our knowledge, the first systematic treatment of this model is carried out by David Aldous \cite{aldous_lecture_2012}, which is followed up by a very recent study presented in \cite{chatterjee_phase_2019}. In both works, the focus is related to the estimation of the so-called mixing times and hence the targeted audiences are mathematicians from the Markov chain mixing time community. In this manuscript, we intend to give a kinetic theory perspective of the model. Indeed, under the large population $N \to \infty$ limit, We can rigorously show that the law of the wealth of a typical agent
(say $X_1$) satisfies the following limit PDE in a weak sense:
\begin{equation}\label{eq:PDE}
\partial_t \rho(t,x) = 2(\rho*\rho)(t,2x) - \rho(t,x).
\end{equation}

\begin{figure}[!htb]
\centering
\includegraphics[scale=0.9]{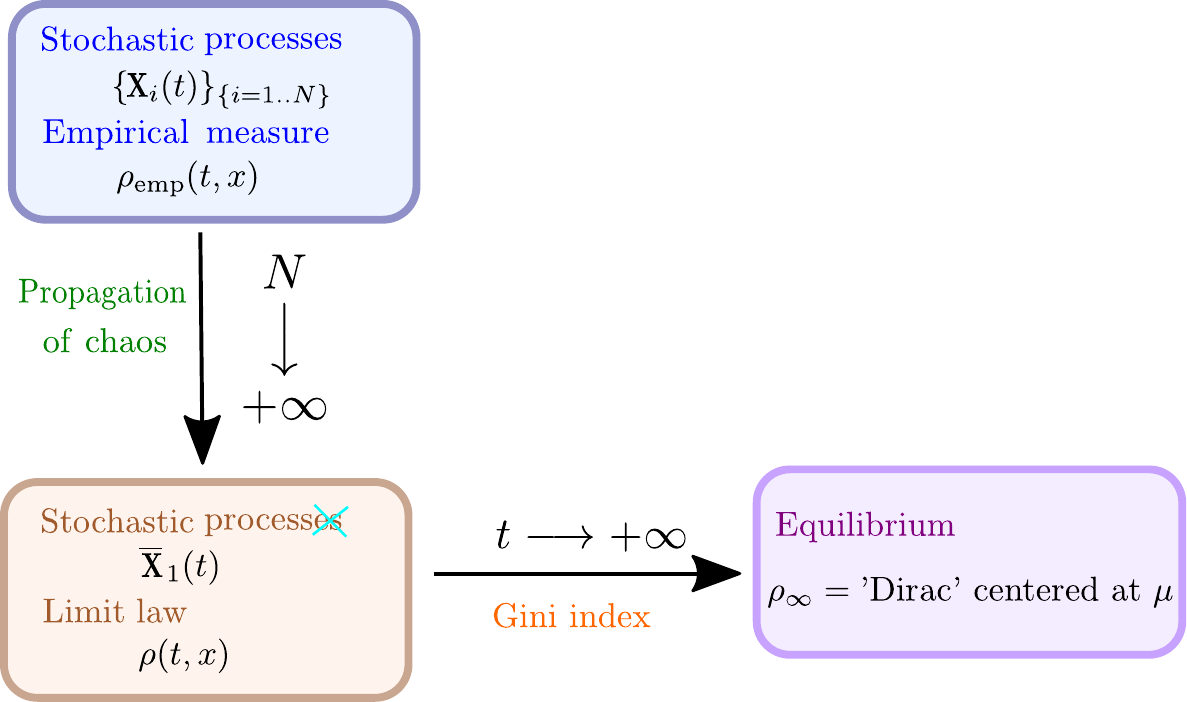}
\caption{Schematic illustration of the general strategy of our treatment of the repeated averaging dynamics, where $\mu$ represents the initial average wealth.}
\label{schema}
\end{figure}

Once the limit PDE is identified from the interacting particle system, the natural next step is to study the problem of convergence to equilibrium of the PDE at hand. In the present work, we demonstrate that the Gini index of $\rho(t)$ converges to $0$ (its minimum value), whence showing that a Dirac distribution centered at the initial average wealth is the equilibrium distribution. Moreover, this model can be served as the first example for which quantitative estimates on the convergence of Gini index can be obtained, which is our primary motivation for writing this paper. An illustration of the general strategy used in this work is shown in Figure \ref{schema}.

Although only a very specific binary exchange model is explored in the present paper, other exchange rules can also be imposed and studied, leading to different models. To name a few, the so-called immediate exchange model introduced in \cite{heinsalu_kinetic_2014} assumes that pairs of agents are randomly and uniformly picked at each random time, and each of the agents transfer a random fraction of its money to the other agents, where these fractions are independent and uniformly distributed in $[0,1]$. The so-called uniform reshuffling model investigated in \cite{dragulescu_statistical_2000,lanchier_rigorous_2018,cao_entropy_2021} suggests that the total amount of money of two randomly and uniformly picked agents possess before interaction is uniformly redistributed among the two agents after interaction. For models with saving propensity and with debts, we refer the readers to \cite{chakraborti_statistical_2000}, \cite{chatterjee_pareto_2004} and \cite{lanchier_rigorous_2018-1}. Also, one can also modulate the rule of picking agents, leading to biased models of money exchange, see for instance \cite{cao_derivation_2021}.

This manuscript is organized as follows: in section \ref{sec:limit_pde}, we briefly discuss the heuristic derivation of the limit equation \eqref{eq:PDE} and give convergence results for the solution of \eqref{eq:PDE} in terms of variance of the distribution as well as the Gini index. Finally, we draw a conclusion in section \ref{sec:conclusion} and present a rigorous treatment of the propagation of chaos phenomenon in Appendix \ref{sec:Appendix}, by applying the martingale-based technique employed in \cite{cao_entropy_2021}.

\section{Convergence to Dirac distribution}
\label{sec:limit_pde}
\setcounter{equation}{0}

We present a heuristic argument behind the derivation of the limit PDE \eqref{eq:PDE} arising from the repeated averaging dynamics in section \ref{subsec:2.1}, we also record a useful stochastic representation behind the PDE \eqref{eq:PDE} on which we will heavily rely. Section \ref{subsec:2.2} is devoted to the exponential decay of the variance of the solution $\rho(t)$ of \eqref{eq:PDE}. In section \ref{subsec:2.3}, we establish a quantitative convergence result on the Gini index of the probability distribution $\rho(t)$, by enforcing the log-concavity property of the initial datum.

\subsection{Formal derivation of the limit PDE}
\label{subsec:2.1}

Introducing $\mathrm{N}^{(i,j)}_t$ independent Poisson processes with intensity $1/N$, the dynamics can be written as:
\begin{equation}
\label{eq:SDE}
\dd X_i(t) = \sum_{j=1..N,j\neq i} \left(\frac{X_i(t-)\!+\!X_j(t-)}{2} - X_i(t-)\right) \dd \mathrm{N}^{(i,j)}_t.
\end{equation}
As the number of players $N$ goes to infinity, one could expect that the processes $X_i(t)$ become independent and of same law. Therefore, the limit dynamics would be of the form:
\begin{equation}
\label{eq:limit_SDE}
\dd \overbar{X}(t) = \left(\frac{\overbar{X}(t-)\!+\!\overbar{Y}(t-)}{2} - \overbar{X}(t-)\right) \dd \overbar{\mathrm{N}}_t,
\end{equation}
where $\overbar{Y}(t)$ is an independent copy of $\overbar{X}(t)$ and $\overbar{\mathrm{N}}_t$ a Poisson process with intensity $1$. Taking a test function $φ$, the weak formulation of the dynamics is given by:
\begin{equation}
\label{eq:limit_SDE_weak}
\dd 𝔼[φ(\overbar{X}(t))] = 𝔼\left[φ\left(\frac{\overbar{X}(t)\!+\!\overbar{Y}(t)}{2}\right) - φ(\overbar{X}(t))\right]\,\dd t.
\end{equation}
In short, the limit dynamics correspond to the jump process:
\begin{equation}
\label{eq:limit_process}
\overbar{X} \leadsto \frac{\overbar{X}\!+\!\overbar{Y}}{2}.
\end{equation}
Let us denote $\rho(t,x)$ the law of the process $\overline{X}(t)$. To derive the evolution equation for $\rho(t,x)$, we need to translate the effect of the jump of $\overbar{X}(t)$ via \eqref{eq:limit_process} onto $\rho(t,x)$.

\begin{lemma}\label{lem1}
Suppose $X$ and $Y$ two independent random variables with probability density $\rho(x)$ supported on $[0,\infty)$. Let $Z = (X+Y)/2$, then the law of $Z$ is given by $Q_+[\rho]$ with:
\begin{equation}
\label{eq:Q_plus}
Q_+[\rho](x) = 2(\rho*\rho)(2x),\quad \forall x\geq 0.
\end{equation}
\end{lemma}

\noindent The proof of this lemma is quite elementary and will be omitted. We can now write the evolution equation for the law of $\overbar{X}(t)$, the density $\rho(t,x)$ satisfies weakly:
\begin{equation}
\label{PDE}
\partial_t \rho(t,x) = \mathcal{L}[\rho](t,x) \qquad \text{ for } t\geq 0 ~~\text{and}~~ x\geq 0
\end{equation}
with
\begin{equation}
\label{eq:L}
\mathcal{L}[\rho](x) := Q_+[\rho](x)-\rho(x) = 2(\rho*\rho)(2x) - \rho(x).
\end{equation}

\begin{remark}
Suppose that $\rho(0,x)$ is a probability density on $[0,\infty)$ with mean $\mu > 0$. It is readily checked that the dynamics \eqref{PDE} preserves the total mass and the mean value. That is, \[\frac{\dd}{\dd t}\int_{\mathbb R_+} \rho(t,x)\, \dd x = 0 ~~\text{and}~~ \frac{\dd}{\dd t}\int_{\mathbb R_+} x\,\rho(t,x)\, \dd x = 0.\] For each test function $φ$, one can show that $\int_{\mathbb R_+} φ(x)\,G[\delta_\mu](dx) = 0$, implying that the Dirac distribution centered at $\mu$ is a equilibrium solution of \eqref{PDE}.
\end{remark}

We now present a stochastic representation of the evolution equation \eqref{PDE}, which is interesting in its own right.

\begin{proposition}\label{stochastic_representation}
Assume that $\rho_t(x) := \rho(t,x)$ is a solution of \eqref{PDE} with initial condition $\rho_0(x)$ being a probability density function supported on $\mathbb{R}_+$ with mean $\mu$. Defining $(X_t)_{t\geq 0}$ to be a $\mathbb R_+$-valued continuous-time pure jump process with jumps of the form
\begin{equation}\label{eq:jump}
\begin{array}{ccc}
X_t & \begin{tikzpicture} \draw [->,decorate,decoration={snake,amplitude=.4mm,segment length=2mm,post length=1mm}]
(0,0) -- (1,0); \node[above,red] at (0.5,0) {\tiny{rate 1}};\end{tikzpicture} & \frac{X_t+Y_t}{2},
\end{array}
\end{equation}
where $Y_t$ is a i.i.d. copy of $X_t$, and the jump occurs according to a Poisson clock running at the unit rate. If $\mathrm{Law}(X_0) = \rho_0$, then $\mathrm{Law}(X_t) = \rho_t$ for all $t\geq 0$.
\end{proposition}

\Proof
Taking $\varphi$ to be an arbitrary but fixed test function, we have
\begin{equation}\label{eq:testfunc}
\frac{\dd }{\dd t} \mathbb E[\varphi(X_t)] = \mathbb E[\varphi((X_t+Y_t)/2)] - \mathbb E[\varphi(X_t)].
\end{equation}
Denoting $\rho(t,x)$ as the probability density function of $X_t$, \eqref{eq:testfunc} can be rewritten as
\[\frac{\dd }{\dd t} \int_{\mathbb{R}_+} \rho(t,x)\varphi(x)\, \dd x = \int_{\mathbb{R}^2_+} \varphi((k+\ell)/2)\rho(k,t)\rho(\ell,t) \,\dd k\, \dd \ell - \int_{\mathbb{R}_+} \rho(t,x)\varphi(x)\, \dd x.\] After a simple change of variables, one arrives at
\begin{equation*}
\frac{\dd }{\dd t} \int_{\mathbb{R}_+} \rho(t,x)\varphi(x)\, \dd x = \int_{\mathbb{R}_+} \left(Q_+[\rho](x,t) - \rho(t,x)\right)\varphi(x)\,\dd x.
\end{equation*}
Thus, $\rho$ has to satisfy $\partial_t \rho = \mathcal{L}[\rho]$ and the proof is completed. \qed

\subsection{Exponential decay of the variance}
\label{subsec:2.2}

Our main goal in this subsection is the proof of the following
\begin{theorem}\label{thm1}
Assume that $\rho(t,x)$ is a classical solution of \eqref{PDE} for each $t >0$, with the initial condition $\rho(0,x)$ being a probability density on $[0,\infty)$ with mean $\mu > 0$ and finite variance. Then the variance of $\rho$ at time $t$, denoted by $\VV(t)$, decays exponentially in time. More specifically, we have $\VV(t) = \VV(0)\,\expo^{-\frac 12t}$.
\end{theorem}

\Proof
Thanks to the conservation of the mean value, we have \[\VV(t) = \int_{\mathbb R_+} x^2\,\rho(t,x)\, dx - \mu^2.\] Thus, we deduce
\begin{align*}
\frac{\dd}{\dd t}\VV(t) &= 2\int_{\mathbb{R}_+} x^2\,(\rho*\rho)(2x)\,\dd x - \int_{\mathbb{R}_+} x^2\,\rho(x)\,\dd x \\
&=\int_{\mathbb{R}_+} 2\,x^2\left(\int_0^{2x} \rho(y)\,\rho(2x-y)\,\dd y\right)\,\dd x - \int_{\mathbb{R}_+} x^2\,\rho(x)\,\dd x \\
&=\int_{y\geq 0}\rho(y)\left(\int_{x\geq y/2} 2\,x^2\,\rho(2x-y)\, \dd x\right)\,\dd y - \int_{\mathbb{R}_+} x^2\,\rho(x)\,\dd x \\
&=\int_{y\geq 0}\rho(y)\left(\int_{z\geq 0} ((y+z)/2)^2\,\rho(z)\, \dd z\right)\,\dd y - \int_{\mathbb{R}_+} x^2\,\rho(x)\,\dd x \\
&= -\frac 12\left(\int_{\mathbb R_+} x^2\,\rho(t,x)\, dx - \mu^2\right) = -\frac 12\VV(t).
\end{align*}
A simple integration yields the advertised conclusion.

\begin{remark}
The proof of Theorem \ref{thm1} can also be carried out from a purely stochastic point of view, by leveraging the stochastic representation of the PDE \eqref{PDE}. Indeed, suppose that $(X_t)_{t\geq 0}$ and $(Y_t)_{t\geq 0}$ are defined as in the statement of Proposition \ref{stochastic_representation}. Then we can calculate
\[\frac{\dd}{\dd t}\VV(t) = \frac{\dd}{\dd t}\Var[X_t] = \Var[(X_t+Y_t)/2] - \Var[X_t] = -\frac 12\Var[X_t] = -\frac 12\VV(t),\] which leads us to the same result.
\end{remark}

% · Var[X_t] = E[|X_t - μ|^2] = W^2_2(ρ(t), δ_μ), this shows the exponential decay of W_2(ρ(t), δ_μ) towards 0 as well.

\subsection{Exponential decay of the Gini index}
\label{subsec:2.3}

The widely used inequality indicator Gini index $G$ measures the inequality in the wealth distribution and ranges from $0$
(no inequality) to $1$ (extreme inequality). We recall the definition of $G$ here for the reader's convenience.

\begin{definition}
Given a probability density function $\rho$ supported on $\mathbb{R}_+$ with mean value $\mu >0$. The Gini index of $\rho$ is given by \[G[\rho] = \frac{1}{2\mu}\iint_{\mathbb{R}^2_+} \rho(x)\,\rho(y)\,|x-y|\,\dd x\,\dd y.\] Alternatively, we can also rewrite \[G[\rho] = \frac{1}{2\mu}\mathbb{E}[|X-Y|],\] in which $X$ and $Y$ are i.i.d. random variables with law $\rho$.
\end{definition}

In econophysics literature, analytical results on Gini index are comparatively rare. In certain models, the Gini index can be shown to converge to $1$, which implies the emergence of the "rich-get-richer" phenomenon and the accentuation of the wealth inequality, see for instance \cite{boghosian_oligarchy_2017,boghosian_h_2015} and references therein. There is also a recently proposed model known as the rich-biased model \cite{cao_derivation_2021}, in which the authors observe a numerical evidence for the convergence of Gini index to its maximum possible value but analytical justification is still absent. As have been indicated earlier, the limit PDE \eqref{PDE} associated with the repeated averaging model can be served as the first example for which quantitative estimates on the behavior of Gini index can be hoped. We start with the following preliminary observation.

\begin{proposition}\label{prop}
Assume that $\rho(t,x)$ is a classical solution of \eqref{PDE} for each $t >0$, with the initial condition $\rho(0,x)$ being a probability density on $[0,\infty)$ with mean $\mu > 0$. Then the Gini index $G[\rho]$ is non-increasing in time. Moreover, we have
\begin{equation}
\label{eq:evolution_of_G}
\begin{aligned}
\frac{\dd}{\dd t}G[\rho] &= -\frac{1}{\mu}\iiint_{\mathbb{R}^3_+} \rho(v)\,\rho(w)\,\rho(y)\left(\frac{|v-y|+|w-y|}{2} - \left|\frac{v+w}{2}-y\right|\right)\,\dd v\,\dd w\,\dd y \\
&\leq 0.
\end{aligned}
\end{equation}
\end{proposition}

\Proof
By symmetry, we have
\begin{align*}
\frac{\dd}{\dd t}G[\rho] &= \frac{1}{\mu}\iint_{\mathbb{R}^2_+} \partial_t\rho(x)\,\rho(y)\,|x-y|\,\dd x\,\dd y \\
&= \frac{1}{\mu}\iint_{\mathbb{R}^2_+} 2\,(\rho*\rho)(2x)\,\rho(y)\,|x-y|\,\dd x\,\dd y - 2\,G[\rho] \\
&= \frac{1}{\mu}\iint_{\mathbb{R}^2_+} 2\left(\int_0^{2x} \rho(z)\,\rho(2x-z)\,\dd z\right)\,\rho(y)\,|x-y|\,\dd x\,\dd y - 2\,G[\rho] \\
&= \frac{1}{\mu}\iiint_{\mathbb{R}^3_+} \rho(v)\,\rho(w)\,\rho(y)\,\left|\frac{v+w}{2} - y\right|\,\dd v\,\dd w\,\dd y - 2\,G[\rho] \\
&= -\frac{1}{\mu}\iiint_{\mathbb{R}^3_+} \rho(v)\,\rho(w)\,\rho(y)\left(\frac{|v-y|+|w-y|}{2} - \left|\frac{v+w}{2}-y\right|\right)\,\dd v\,\dd w\,\dd y,
\end{align*}
whence the proof is finished. \qed

\begin{remark}
In light of the previous remark and stochastic representation of the PDE \eqref{PDE}. We can also provide an alternative proof of Proposition \ref{prop}. Indeed, suppose that $(X_t)_{t\geq 0}$ and $(Y_t)_{t\geq 0}$ are defined as in the statement of Proposition \ref{stochastic_representation}. Then we can compute
\begin{align*}
\frac{\dd}{\dd t}G[\rho] = \frac{1}{2\mu}\frac{\dd}{\dd t}\mathbb{E}[|X_t - Y_t|] &= \frac{1}{\mu}\mathbb{E}[|(X_t+Z_t)/2 - Y_t|] - \frac{1}{\mu}\mathbb{E}[|X_t-Y_t|] \\
&= -\frac{1}{\mu}\left(\mathbb{E}[|X_t-Y_t|] - \mathbb{E}[|(X_t+Z_t)/2 - Y_t|]\right) \leq 0,
\end{align*}
in which $Z_t$ is a fresh i.i.d. copy of $X_t$ (independent of $Y_t$ as well). This coincides with \eqref{eq:evolution_of_G}
\end{remark}

At this point, we may expect to bound $G[\rho]$ in terms of $-\frac{\dd}{\dd t} G[\rho]$ in order to extract some information on the rate of decay of $G$. But unfortunately, inequalities of the form $-\frac{\dd}{\dd t} G[\rho] \geq c\cdot G[\rho]$ can not be always fulfilled. For example, if we take $\rho = \frac 12\delta_0 + \frac 12\delta_{2\mu}$, then one can check that $\frac{\dd}{\dd t} G[\rho] = 0$, whereas $G[\rho] = \frac 12 > 0$. However, not all hope is lost. Indeed, if we restrict the initial data $\rho(0,x)$ to be log-concave, we can prove the following

\begin{theorem}\label{main_thm}
Assume that $\rho(t,x)$ is a classical solution of \eqref{PDE} for each $t >0$, with the initial condition $\rho(0,x)$ being a log-concave probability density on $[0,\infty)$ with mean $\mu > 0$. Then the Gini index $G[\rho]$ converges to $0$ exponentially fast in time. Moreover, we have
\begin{equation}
\label{eq:decay_of_G}
G[\rho(t)] \leq G[\rho(0)]\expo^{-\frac{t}{14434}}.
\end{equation}
\end{theorem}

To facilitate the proof of Theorem \ref{main_thm}, we need the following

\begin{lemma}\label{log-concavity}
Assume that $\rho(t,x)$ is a classical solution of \eqref{PDE} for each $t >0$, with the initial condition $\rho(0)$ being a log-concave probability density on $[0,\infty)$ with mean $\mu > 0$. Then $\rho(t)$ is again log-concave for each $t>0$.
\end{lemma}

\Proof The proof is an immediate consequence of the stochastic representation of \eqref{PDE}, together with the elementary fact that log-concavity is preserved by convolution. \qed

\begin{remark}
Preservation of log-concavity can also be established for other PDEs, although the proofs are usually quite involved. For instance, it is well-known that evolution under the one-dimensional heat equation preserves the log-concavity of the initial datum \cite{brascamp_extensions_2002}.
\end{remark}

\noindent{\bf Proof} of Theorem \ref{main_thm} \quad For notational simplicity, we write
\begin{equation}\label{notations}
G:=G[\rho] ~~ \text{and}~~ H:= -\frac{\dd}{\dd t} G[\rho].
\end{equation}

In fact, we will not need the restriction that the support of the distribution $\rho$ is $[0,\infty)$.

By approximation, without loss of generality, we may assume that $\rho(x)>0$ for all real $x$. For example, one may approximate $\rho$ by its convolution $p*\phi$ with the density $\phi$ of a centered normal distribution with an arbitrarily small variance. Then $\rho*g>0$ on $\mathbb R$ and $\rho*\phi$ is arbitrarily close to $\rho$ and log-concave, thanks to the preservation of log-concavity by convolution.

\begin{figure}[ht]
\centering
\includegraphics[scale=0.8]{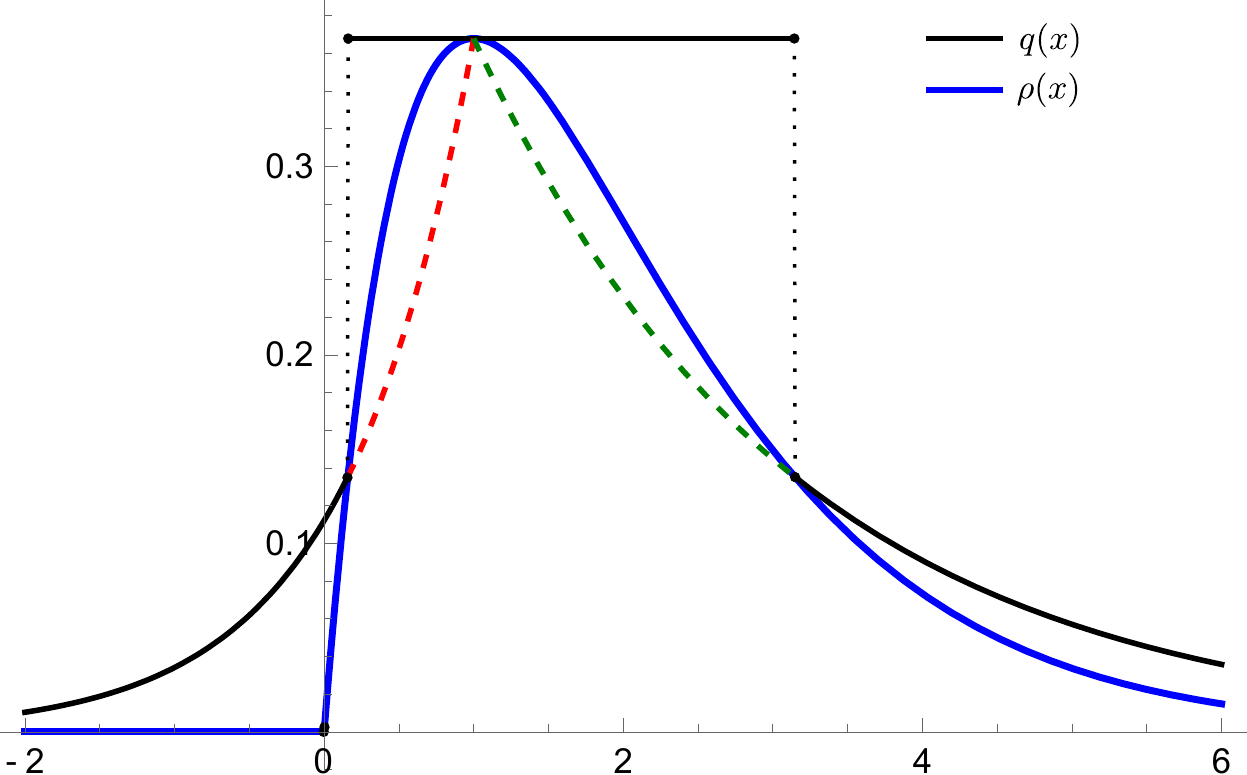}%[width=.35\textwidth]
\caption{For $\rho(x)=x\,\expo^{-x}\,\mathbbm{1}_{\{x>0\}}$, here are the graphs $\{(x,\rho(x)) \mid -2\leq x\leq 6\}$ (blue), $\{(x,q(x))\mid -2\leq x\leq 6\}$ (black), $\{(x,q_1(x))\mid a\leq x\leq c\}$ (dashed red), and $\{(x,q_2(x))\mid c\leq x\leq b\}$ (dashed green). For this particular $\rho$, we have $c=1$, $a=-W_0\left(-1/\expo^2\right) \approx 0.1586$, and $b=-W_{-1}\left(-1/\expo^2\right) \approx 3.1461$, where $W_j$ is the $j$th branch of the Lambert $W$ function \cite{Lambert_Observationes_1758}.}
\label{fig:q}
\end{figure}

As $\rho$ is a log-concave density, $\rho$ is continuous and attains its maximum value, say $\rho_* (>0)$, at some point $c\in \mathbb R$, so that $\rho_*=\rho(c)\geq \rho(x)$ for all real $x$. Moreover, again because $\rho$ is log-concave, there exist (unique) real $a$ and $b$ such that
\begin{equation*}
a<c<b ~~\text{and}~~ \rho(a)=\rho(b)=\rho_*/\expo.	
\end{equation*}
We define \begin{equation*}
q(x):=
\left\{
\begin{aligned}
&q_1(x):=\rho_*\exp\left\{-\frac{x-c}{a-c}\right\} &~\text{ if }x<a, \\
&\rho_* &~\text{ if }a\leq x<b, \\
&q_2(x):=\rho_*\exp\left\{-\frac{x-c}{b-c}\right\} &~\text{ if }x\geq  b.
\end{aligned}
\right.
\end{equation*}
Thanks to the log-concavity of $\rho$ again, we have $\rho(x) \leq q(x)$.  We refer to Figure \ref{fig:q} for an illustration.

By shifting, we may assume with of loss of generality that $a=0$. Thus,
\begin{align*}
G&\leq \iint\limits_{\mathbb{R}^2} q(x)\,q(y)\,|x-y|\,\dd x\,\dd y \\
&=\rho_*^2\frac{\expo^2 b^3+9 \expo b^3+3 b^3+3 b^2 c-12 \expo b^2 c-3 b c^2+12 \expo b c^2}{3 \expo^2} \\
&\leq \rho_*^2\frac{\left(1+3 \expo + \expo^2/3\right) b^3}{\expo^2},
\end{align*}
since $0<c<b$.
Moreover, again by the log-concavity of $p$, we have $\rho \geq \rho_*/\expo$ on the interval $[a,b]=[0,b]$, so that $1=\int_{\mathbb R} \rho \geq \int_0^b \rho_*/\expo = b\,\rho_*/\expo$, whence $\rho_*\leq \expo/b$ and
\begin{equation}\label{eq:upper_G}
G\leq (1+3\expo +\expo^2/3)\,b.
\end{equation}

On the other hand, because $\rho \geq \rho_*/\expo$ on the interval $[a,b]=[0,b]$ and the integrand in the definition of $H$ is non-negative, we have
\begin{align*}
H&\geq\left(\frac{p_*}{\expo}\right)^3\iiint\limits_{[0,b]^3} \left(\frac{|x-z|+|y-z|}2
-\left|\frac{x-z+y-z}2\right|\right)\,\dd x\,\dd y\,\dd z \\
&=\left(\frac{\rho_*}{\expo}\right)^3\frac{b^4}{24}.
\end{align*}
Also, $1=\int_{\mathbb R} \rho\leq \int_{\mathbb R} q=\rho_*b(1+1/\expo)$, so that $\rho_*\geq 1/(b(1+1/\expo))$ and hence
\begin{equation}\label{eq:lower_H}
H\geq \left(\frac1{(\expo +1)b}\right)^3\frac{b^4}{24}=\frac b{24(\expo + 1)^3}.
\end{equation}

\noindent Comparing \eqref{eq:upper_G} and \eqref{eq:lower_H}, we deduce
\begin{equation*}
H\geq \frac G{24(\expo + 1)^3(1+3 \expo + \expo^2/3)}\geq \frac{G}{14334},
\end{equation*}
as claimed. \qed

\begin{figure}[!htb]
\centering
\includegraphics[scale = 0.8]{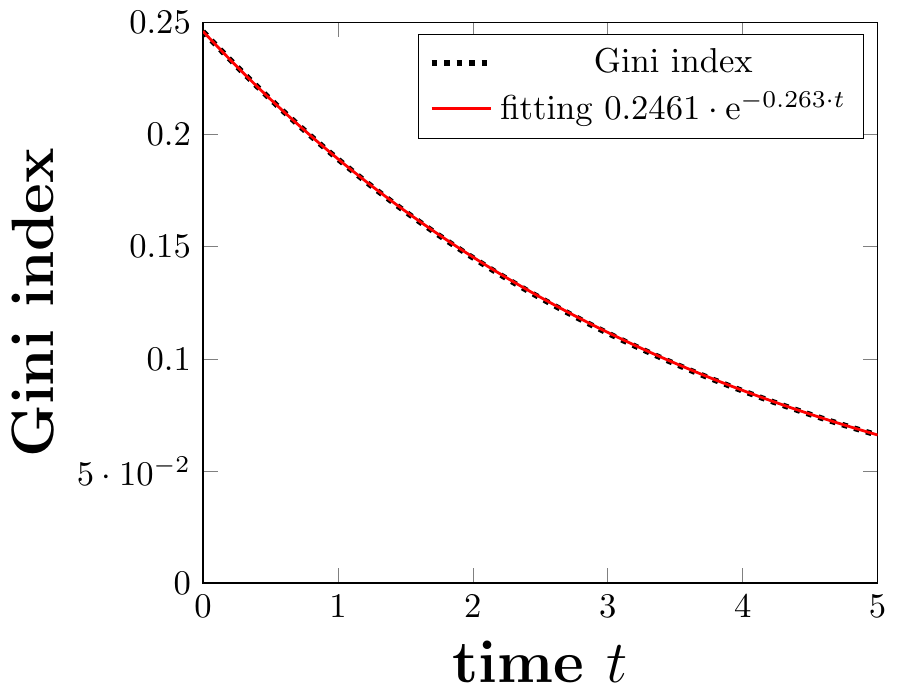}
\caption{Evolution of the Gini index of $\rho(t)$ (the solution of \eqref{PDE}) for $0\leq t \leq 5$, with the initial datum being a gamma probability density with shape parameter $\mu = 5$ and rate parameter equal to unity, i.e., $\rho(0,x) = \mathbbm{1}_{[0,\infty)}(x)\cdot x^{\mu - 1}\,\expo^{-x} / \Gamma(\mu)$. The black dotted line and the red smooth curve represent the Gini index and its fitting curve, respectively. We also remark that in this experiment these two curves are almost indistinguishable.}
\label{Gini_evolution}
\end{figure}

Finally, we provide a numerical experiment in order to corroborate the relaxation of the Gini index guaranteed by Theorem \ref{main_thm}, see Figure \ref{Gini_evolution}. For the initial condition, we use a gamma probability density with shape parameter $\mu = 5$ and rate parameter equal to unity, i.e., $\rho(0,x) = \mathbbm{1}_{[0,\infty)}(x)\cdot x^{\mu - 1}\,\expo^{-x} / \Gamma(\mu)$. The standard forward Euler scheme (with the time step-size $\Delta t = 0.05$ and the space step-size $\Delta x = 0.01$) is enforced for the numerical solution of \eqref{PDE}. Note that the Gini index of our choice of $\rho(0,x)$ has a nice closed expression $G[\rho(0)] = \frac{2^{1-2\mu}\,\Gamma(2\mu)}{\mu\,(\Gamma(\mu))^2}$, which reduces (approximately) to 0.2461 for $\mu = 5$.

\section{Conclusion}
\label{sec:conclusion}

In this manuscript, we have investigated the repeated averaging dynamics for money exchange originated from econophysics. Because of the model simplicity in its appearance, there is a comparative lack of mathematical literature that is purely dedicated to this model, although this model is a special case of the general dynamics studied in \cite{matthes_steady_2008}. We presented a propagation of chaos result, which links the stochastic $N$ particle system to a deterministic nonlinear evolution equation. Although certain convergence results of the Gini index are obtained for other econophysics models, we emphasize that no quantitative estimates on the long time behavior of Gini index are available in the current literature (at least to our best knowledge). Thus, this toy model may serve as a starting point for more systematic, quantitative investigation of the large-time asymptotic of Gini index arising from other models. As a open conjecture, we speculate that the constant $1/14434$ appearing in the statement of Theorem \ref{main_thm} might be tremendously improved.

It would also be interesting to investigate the behavior of the Gini index for the stochastic agent-based model where the number of agents $N$ is arbitrary but fixed. We believe that it would be relatively simple (in this setting) to demonstrate the convergence of the Gini index towards zero, but the difficulty arises when we want to obtain an explicit rate of the aforementioned convergence.\\

\noindent {\bf Acknowledgement~} It is a great pleasure to thank Iosif Pinelis for his answer to a question of myself on MathOverflow \cite{pinelis_possibility_2021}, where the essential piece needed for a complete proof of Theorem \ref{main_thm} is presented. I would also want to express my gratitude to my Ph.D advisor Sebastien Motsch for his careful proofreading of the manuscript.

\section{Appendix}

\subsection{Propagation of chaos}
\label{sec:Appendix}
In the last part of the manuscript, we sketch the proof of the so-called propagation of chaos \cite{sznitman_topics_1991,chaintron_propagation_2021}, relying on a  martingale-based technique developed in \cite{merle_cutoff_2019}. We emphasize that the proof presented here only a slight modification of the Theorem 6 in \cite{cao_entropy_2021}.

We equip the space $\mathcal{P}(\mathbb{R}_+)$ with the Wasserstein distance with exponent $1$, which is defined via \[W_1(\mu,\nu) = \sup\limits_{\|\nabla \varphi\|_\infty \leq 1} \langle \mu-\nu, \varphi\rangle\] for $\mu,\nu \in \mathcal{P}(\mathbb{R}_+)$. The propagation of chaos result is summarized in the following

\begin{theorem}\label{PoC}
Denote the empirical distribution of the repeated averaging $N$ particle system \eqref{eq:repeated_averaging} at time $t$ as \[\rho_{\emp}(t):= \frac{1}{N}\,\sum_{i=1}^N \delta_{X_i(t)},\] and let $\rho(t)$ be the solution of \eqref{PDE} with initial data $\rho(0)$. If \begin{equation}\label{assump}
\mathbb{E}[W_1(\rho_{\mathrm{\emp}}(0),\rho(0))] \xrightarrow[]{} 0 ~\text{as}~ N\to \infty,
\end{equation}
then we have that\[\mathbb{E}[W_1(\rho_{\mathrm{\emp}}(t),q(t))] \xrightarrow[]{} 0 ~\text{as}~ N\to \infty,\] holding for all $0\leq t \leq T$ with any prefixed $T >0$.
\end{theorem}

\Proof We recall that the map $Q_+[\cdot] \colon \mathcal{P}(\mathbb{R}_+) \to \mathcal{P}(\mathbb{R}_+)$ is defined via
\[Q_+[\rho](x) = 2(\rho*\rho)(2x),\quad \forall x\geq 0.\] Assume that a classical solution $\rho(t,x)$ of
\begin{equation}\label{APDE}
\rho(t,x) = \rho(0,x) + \int_0^t \mathcal{L}[\rho](s,x)\,\dd s
\end{equation}
exists for $0\leq t < \infty$, where $\mathcal{L} = Q_+ - \mathrm{Id}$ and $\rho(0,x)$ is a probability density function whose support is contained in $\mathbb{R}_+$. The map $Q_+$ is Lipschitz continuous in the sense that
\begin{equation}\label{LipG}
W_1(Q_+[f],Q_+[g]) \leq W_1(f,g)
\end{equation}
for any $f,g \in \mathcal{P}(\mathbb{R}_+)$. Indeed, we have \[W_1(Q_+[f],Q_+[g]) = \sup\limits_{\|\nabla \varphi\|_\infty \leq 1} \mathbb{E}\left[\varphi((X_1+Y_1)/2) - \varphi((X_2+Y_2)/2)\right],\] where $X_1, Y_1$ are i.i.d with law $f$, $X_2, Y_2$ are i.i.d with law $g$. By Lipschitz continuity of the test function $\varphi$, we obtain
\[W_1(Q_+[f],Q_+[g]) \leq \mathbb{E}[|X_1-X_2|].\] We now recall an alternative formulation of $W_1(f,g)$, given by \[W_1(f,g) = \inf\left\{\mathbb E[|X-Y|];~\mathrm{Law}(X)=f,~\mathrm{Law}(Y)=g\right\},\] so in particular, we may take a coupling of $X_1$ and $X_2$ so that $W_1(f,g) = \mathbb{E}[|X_1-X_2|]$. Assembling these pieces together, we arrive at \eqref{LipG}.  More generally, suppose we have two random probability measures $f$ and $g$ with bounded second moment, taking expectation on both sides of \eqref{LipG} gives rise to
\begin{equation}\label{newLipG}
\mathbb{E} \left[\sup_{\|\nabla\varphi\|_{\infty}\leq 1} \int \varphi(x)\,(Q_+[f]-Q_+[g])\right]
\leq \mathbb{E} \left[\sup_{\|\nabla\varphi\|_{\infty}\leq 1} \int \varphi(x)\,(f(\dd x)-g(\dd x))\right].
\end{equation}
We now observe that the empirical measure is a compound jump process: Define $N_t$ a homogeneous Poisson process with constant intensity $\lambda=(N-1)/2$. Given $\tau_1,\ldots,\tau_k$ the times when $N_t$ jumps, we take the $Y_{\tau_k}$ independent: At each $\tau_k$, with uniform probability $\frac{2}{N\,(N-1)}$ we choose a pair $i< j$ and take
\[\begin{split}
Y_{\tau_k}=&\frac{1}{N}\,\Big(2\delta(x-(X_i(\tau_k-)+X_j(\tau_k-)/2))
\\ &-\delta(x-X_i(\tau_k-))-\delta(x-X_j(\tau_k-))\Big).
\end{split}
\]
We immediately note that
\begin{equation}
\begin{split}
\lambda\,\mathbb{E}[Y_t] &= \frac{1}{N^2}\,\sum_{i<j} \mathbb{E}\Big[2\delta(x- (X_i(t-)+X_j(t-)/2))\\&\ -\delta(x-X_i(t-))-\delta(x-X_j(t-))\Big].
\end{split}\label{expectYt}
\end{equation}

We now show that the empirical measure of the stochastic system satisfies an approximate version of \eqref{APDE}. Fix a deterministic test function $\varphi$ with $\|\nabla \varphi\|_\infty \leq 1$, and consider the time evolution of $\langle \rho_{\emp},\varphi \rangle$ where for some probability measure $\nu$, we denote by the duality bracket $\langle \nu,\varphi \rangle =\int \varphi\,\dd \nu$. Then
\[\dd \mathbb{E}[\langle \rho_{\emp},\varphi \rangle]=\dd \mathbb{E}\left[\langle Y_t\,\dd N_t,\varphi \rangle \right] =\lambda\,\langle\mathbb{E}[Y_t],\varphi \rangle\,\dd t.\]
Therefore, thanks to \eqref{expectYt},
\[
\begin{aligned}
\dd \mathbb{E}[\langle \rho_{\emp},\varphi \rangle] &= \frac{1}{N^2}\,\sum_{i < j} \mathbb{E}\left[2\varphi\left((X_i+X_j)/2\right) + - \varphi(X_i) - \varphi(X_j)\right] \dd t \\
&=\frac{1}{N^2}\,\sum_{i,j=1\ldots N,i\neq j} \mathbb{E}\left[\varphi\left((X_i+X_j)/2\right)-\varphi(X_i)\right] \dd t\\
&=\frac{1}{N^2}\,\sum_{i,j=1}^N \mathbb{E}\left[\varphi\left((X_i+X_j)/2\right)-\varphi(X_i)\right] \dd t,
\end{aligned}
\]
where all $X_i,\,X_j$ are taken at time $t-$. On the other hand, we may calculate
\[\langle Q_+[\rho_{\emp}], \varphi \rangle = \int \varphi(x)\,2\,\frac{1}{N^2}\,\sum_{i,j=1}^N \delta_{X_i+X_j}(2x)\,\dd x = \frac{1}{N^2}\,\sum_{i,j=1}^N \varphi\left((X_i+X_j)/2\right).\]
Therefore
\begin{equation}
\dd \mathbb{E}[\langle \rho_{\emp},\varphi \rangle]=\mathbb{E}\left[\langle \mathcal{L}[\rho_{\emp}], \varphi \rangle\right]\,\dd t.\label{weakformulation}
\end{equation}
By Dynkin's formula, the compensated process
\begin{equation}\label{mart}
M_\varphi(t):= \langle \rho_{\emp}(t),\varphi \rangle - \langle \rho_{\emp}(0),\varphi \rangle - \int_0^t \mathbb{E}[\langle \mathcal{L}[\rho_{\emp}(s)], \varphi \rangle] \, \dd s
\end{equation}
is a martingale. Furthermore, comparing with \eqref{APDE}, we easily obtain that
\[\begin{aligned}
\langle \rho_{\emp}(t) - \rho(t),\varphi \rangle &= M_\varphi(t) + \langle \rho_{\emp}(0) - \rho(0),\varphi \rangle \\
&\quad + \mathbb{E}\,\int_0^t \langle \mathcal{L}[\rho_{\emp}(s)] - \mathcal{L}[\rho(s)], \varphi \rangle \,\dd s.
\end{aligned}\]
Taking the supremum over $\varphi$, we therefore have that
\[\begin{aligned}
&\mathbb{E}\sup_{\|\nabla\varphi\|_{\infty}\leq 1} \langle \rho_{\emp}(t) - \rho(t),\varphi \rangle \leq \mathbb{E}\,\sup_{\|\nabla\varphi\|_{\infty}\leq 1} (|M_\varphi(t)| + \langle \rho_{\emp}(0) - \rho(0),\varphi \rangle) \\
&\quad +\int_0^t \mathbb{E}\,\sup_{\|\nabla\varphi\|_{\infty}\leq 1} \langle \mathcal{L}[\rho_{\emp}(s)] - \mathcal{L}[\rho(s)], \varphi \rangle \,\dd s.
\end{aligned}\]
By the definition of the $W_1$ distance, we deduce from \eqref{newLipG}  that
\[\begin{split}
&\mathbb{E}\,W_1(\rho_{\emp}(t),q(t)) \leq \eta(t) + 2\,\int_0^t \mathbb{E}\, W_1(\rho_{\emp}(t),q(t))\, \dd s,
\end{split}\]
in which we have set
\begin{equation}\label{eta}
\eta(t) :=  \mathbb{E}\,\sup\limits_{\|\nabla \varphi\|_\infty \leq 1}\,|M_\varphi(t)| + \mathbb{E}\,W_1(\rho_{\emp}(0),q(0)).
\end{equation}
Thus, Gronwall's inequality gives rise to
\begin{equation}\label{almostfinish}
\mathbb{E}\,W_1(\rho_{\emp}(t),\rho(t))\leq \left(\sup\limits_{t\in [0,T]} \eta(t)\right)\expo^{2\,T}.
\end{equation}
In order to establish propagation of chaos for $t \leq T$, it therefore suffices to show that
\begin{equation}\label{convprob}
\sup\limits_{t\in [0,T]} \eta(t) \xrightarrow[N \to \infty]{\mathbb{P}} 0.
\end{equation}
To prove \eqref{convprob}, we treat each term appearing in the definition of $\eta(t)$ separately. The second term in \eqref{eta} approaches to 0 as $N \to \infty$ by our assumption. The treatment of the first term is more delicate, but can be carried out in a similar fashion as the proof of Theorem 6 in \cite{cao_entropy_2021}. In the end, we obtain estimates of the form
\[\mathbb{E}\,\left[\sup_{\|\nabla\varphi\|_\infty\leq 1}\, \left|M_\varphi(t)\right| \right]\leq C\,\frac{t^\theta}{N^\theta}\] for some $\theta >0$, which allows to finish the proof of \eqref{convprob}. \qed

\end{document}